\documentclass[12pt]{article}
\usepackage{amsfonts,amsmath,amssymb}
\usepackage{hyperref}
\long\def\nodo#1{}
\def\id{\mathrm{id}}
\def\op{\mathrm{op}}
\def\genfd{\mathbf{k}}
\def\FF{F}

\def\nxpoint{\refstepcounter{subsection}\makepoint{\thesubsection}}
\def\nxsubpoint{\refstepcounter{subsubsection}%
  \makepoint{\thesubsubsection}}
\def\refpoint#1{{\rm\textbf{\ref{#1}}}}
\def\makepoint#1{\medbreak\noindent{\bf #1. }}

\begin{document}
\begin{center}
  {\Large Drinfeld-Xu bialgebroid 2-cocycles twist the antipode
  }
  \\
 \vskip .2in
  {\sc Zoran \v{S}koda}
\end{center}
\begin{abstract}
  Ping Xu generalized Drinfeld 2-cocycles from bialgebras to associative bialgebroids over noncommutative base algebras. Any counital Drinfeld--Xu 2-cocycle twists the base algebra of the bialgebroid and a comultiplication on the total algebra, obtaining a new, twisted bialgebroid. Antipodes for bialgebroids have been considered, but finding a general way to twist the antipode, which is straightforward in the Hopf algebra case, appeared somewhat elusive. In this article, we prove that if an invertible antipode $S$ for the original bialgebroid exists, and another expression $V_\FF$ depending on the 2-cocycle $\FF$ is invertible, then the expected conjugation formula $S_\FF(-) = V_\FF^{-1} S(-) V_\FF$ indeed produces an invertible antipode $S_\FF$ for the twisted bialgebroid.
\end{abstract}

\section{Introduction}

\nxpoint Given an associative bialgebroid $H$ over a noncommutative base $R$, 
Ping Xu~\cite{xu} has introduced a notion of a (nonabelian) 2-cocycle (which he called a ``twistor'') $\FF\in H\otimes H$ and a ``twisting'' procedure where to an $R$-bialgebroid $H$ with a (counital) 2-cocycle $\FF$ one attaches a new bialgebroid $H_\FF$ over a new base $R_\FF$. This is a rather subtle generalization of Drinfeld's twisting of bialgebras by bialgebra 2-cocycles (``Drinfeld twists''). If a bialgebra $H$ is actually a Hopf algebra with antipode $S$, then the twisted bialgebra $H_\FF$ is also Hopf with antipode
$S_\FF: h\mapsto V_\FF^{-1} (Sh) V_\FF$ where, in our conventions for 2-cocycles,
$V_\FF = (S F^{1}) F^{2}$ in Sweedler-like notation for 2-cocycle without summation sign $\FF = F^{1}\otimes F^{2} := \sum_j F^{1}_i\otimes F^{2}_i$. The proof that $S_\FF$ is indeed an antipode is very simple, in part because one can easily show that $V_\FF$ is invertible and write down an explicit formula for $V_\FF^{-1}$ to use it in calculation. For Drinfeld--Xu bialgebroid 2-cocycles, twisting the antipodes in general is an open problem for over 20 years. In this article, I present the first general proof that essentially the same formula for the twisted antipode as in the Hopf algebra case indeed works for the B\"ohm--Szlach\'anyi Hopf algebroids with invertible antipode provided $V_\FF$ is invertible. The proof is not only calculationally nontrivial, but it requires many checks, some involved, that the steps and the intermediate expressions in the calculations are well-defined. 

\nxpoint Our conventions for bialgebroids are as follows. The ground field is $\genfd$. If $R = (R,\mu_R)$ is a (unital associative) $k$-algebra,
then a {\bf left $R$-bialgebroid}
(\cite{bohmnew,bohmHbk,bohmSzlach,BrzMilitaru,lu}) is
a 5-tuple $(H,\alpha,\beta,\Delta,\epsilon)$
where $H = (H,\mu)$ is a $k$-algebra,
source $\alpha:R\to H$ and target $\beta:R^{op}\to H$ are algebra maps
whose images commute inducing an $R$-bimodule structure on $H$ by
$r.h.s = \alpha(r)\beta(s)h$ for $h\in H,r,s\in R$; the comultiplication
$\Delta:H\to H\otimes_R H$, $h\mapsto h_{(1)}\otimes h_{(2)}
= \sum_i h_{(1)i}\otimes h_{(2)i}$,
is an $R$-bimodule map (i.e. $\Delta(\alpha(a)\beta(b)h)=\alpha(a)h_{(1)}\otimes\beta(b)h_{(2)}$) which makes $(H,\Delta)$
a coassociative comonoid in the bicategory of $R$-bimodules
with counit $\epsilon:H\to R$.
Since the monoidal unit coherences
$R\otimes_R H\cong R\cong H\otimes_R R$ are
realized by bimodule isomorphisms $r\otimes_R h\mapsto \alpha(r)h$
and $h\otimes_R r\mapsto \beta(r)h$, 
the counit axioms read (in Sweedler notation)
$\alpha(\epsilon(h_{(1)}))h_{(2)} = h = \beta(\epsilon(h_{(2)}))h_{(1)}$ for
all $h\in H$ and the fact that $\epsilon:H\to R$ is an $R$-bimodule map by
$\epsilon(\alpha(r)h)=r\epsilon(h)$ and $\epsilon(\beta(r)h) = \epsilon(h)r$
and in particular $\epsilon\circ\alpha=\id_R = \epsilon\circ\beta$.
Require that  $\epsilon(g h) = \epsilon(g\alpha(\epsilon(h)))$ for all $g,h\in H$. Since the kernel 
$$
I_R=\operatorname{Ker}(H\otimes_\genfd H\to H\otimes_R H) =
\operatorname{Span}_\genfd\{\beta(r)\otimes 1-1\otimes\alpha(r)| r\in R\}
$$
is only a left ideal in general, $H\otimes_R H$ is not necessarily an algebra.
For this reason, a subtle compatibility of $\Delta$
with the multiplication $\mu:H\otimes_k H\to H$ is demanded~\cite{bohmHbk},
namely $\Delta(h)I_R\subset I_R$ (the left hand side interpreted in
terms of factorwise multiplication in $H\otimes_\genfd H$).
This condition ensures that the rule
$h\otimes(k\otimes_R l)\mapsto (h_{(1)}k)\otimes_R(h_{(2)}l)$ defines a well-defined map $H\otimes_k(H\otimes_R H)\to H\otimes_R H$,
which is furthermore required to be a unital action. 

\nxsubpoint Define a $\genfd$-linear map $\triangleright:H\otimes R\to R$,
$h\otimes r\mapsto h\triangleright r := \epsilon(h\alpha(r))$.
Since $\epsilon\circ\beta=\id_R=\epsilon\circ\alpha$, 
it follows that $h\triangleright r =
h\triangleright\epsilon(\beta(r)) = \epsilon(h\alpha(\epsilon(\beta(r)))) =
\epsilon(h\beta(r))$ and $\alpha(r)\triangleright 1 = r = \beta(r)\triangleright 1$, for $r\in R$ and $h\in H$.
Axiom $\epsilon(g h) = \epsilon(g\alpha(\epsilon(h)))$ just says
$(g h)\triangleright 1_R = g\triangleright(h\triangleright 1_R)$.
Using this property and the associativity of the multiplication in $H$,
the action axiom for $\triangleright$ follows,
$(g h)\triangleright r =
\epsilon(g h\alpha(r)) = \epsilon(g\alpha(\epsilon(h\alpha(r)))) =
g\triangleright(h\triangleright r)$. The unitality of action $1_H\triangleright r = r$ reads $\epsilon\circ\alpha =\id_R$. Conversely, the action axiom for $\triangleright$, the relation between $\triangleright$ and $\epsilon$ and equalities $\epsilon\circ\beta=\id_R=\epsilon\circ\alpha$ imply the
$R$-bilinearity of $\epsilon$, 
by $\epsilon(\alpha(r)h) = \epsilon(\alpha(r)\alpha(\epsilon(h)))
= \epsilon(\alpha(r\epsilon(h))) = r\epsilon(h)$ and
$\epsilon(\beta(r)h) = \epsilon(\beta(r)\beta(\epsilon(h))) =
\epsilon(\beta(\epsilon(h)r)) = \epsilon(h)r$. In particular,
\begin{equation}\label{eq:htrrbilin}
  (\alpha(a)\beta(b) h)\triangleright r = a (h\triangleright r) b.
\end{equation}
  \nxpoint {\bf Lemma.} For all $h\in H$ and $r\in R$
\begin{equation}\label{eq:halphar}\begin{array}{lcl}
h\alpha(r) = \alpha(h_{(1)}\triangleright r)h_{(2)}
\\
h\beta(r) = \beta(h_{(2)}\triangleright r)h_{(1)}
\end{array}\end{equation}

{\it Proof.} $\Delta$ is an $R$-bimodule map, hence
$\Delta(\alpha(r)) = \alpha(r)\otimes_R 1_H$ and
$\Delta(\beta(r)) = 1_H\otimes_R\beta(r)$ for all $r\in R$.
It follows that $\Delta(h\alpha(r)) = \Delta(h)(\alpha(r)\otimes_R 1) =
h_{(1)}\alpha(r)\otimes_R h_{(2)}$ and $\Delta(h\beta(r)) = h_{(1)}\otimes_R h_{(2)}\beta(r)$. The counit axiom then gives
$h\alpha(r) = \alpha(\epsilon(h_{(1)}\alpha(r)))h_{(2)} = \alpha(h_{(1)}\triangleright r)h_{(2)}$ and likewise $h\beta(r) = \beta(\epsilon(h_{(2)}\beta(r)))h_{(1)}
= \beta(h_{(2)}\triangleright r)h_{(1)}$.

\nxpoint Usually one defines a Hopf algebroid as
a bialgebroid with an antipode (though some more general variants exist).
Several nonequivalent versions of the axioms for the defining
antipode of a Hopf algebroid are used in the literature
(see e.g.~\cite{bohmHbk,lu,twosha}).
We adopt Definition~4.1 from~\cite{bohmSzlach}.

A (left) {\bf Hopf $R$-algebroid with an invertible antipode}
is a left $R$-bialgebroid $(H,\alpha,\beta,\Delta,\epsilon)$ with an antipode
map $S: H\to H$ which is invertible antihomomorphism of algebras satisfying
for all $h\in H$ the conditions
\begin{equation}\label{eq:antSba}
S\circ\beta = \alpha
\end{equation}
\begin{equation}\label{eq:antShS}
(S h_{(1)})_{(1)} h_{(2)}\otimes_R(S h_{(1)})_{(2)} = 1_H\otimes_R S h.
\end{equation}
\begin{equation}\label{eq:antSmSmh}
(S^{-1} h_{(2)})_{(1)}\otimes_R(S^{-1} h_{(2)})_{(2)}h_{(1)} = S^{-1} h\otimes_R 1_H
\end{equation}
Axiom~(\ref{eq:antSba}) implies that $(S\otimes \id)I_R \subset I_R$ and
$\mu(S\otimes\id)I_R = 0$.

\nxpoint \label{antihom}
(Proposition 4.4 in~\cite{bohmHbk}) The antipode map is automatically an antihomomorphism of corings, 
$(S h)_{(1)} \otimes_R (S h)_{(2)} = S(h_{(2)})\otimes_R S(h_{(1)})$.

\nxpoint Let $H$ be a left associative $R$-bialgebroid. An element
$\FF\in H\otimes_R H$ is called a 2-{\bf cocycle} if
\begin{equation}\label{eq:2coc}
[(\Delta\otimes_R \mathrm{id})(\FF)](\FF\otimes_\genfd 1)
=
[(\mathrm{id}\otimes_R\Delta)(\FF)](1\otimes_\genfd\FF)  
\end{equation}
in $H\otimes_R H\otimes_R H$.
A 2-cocycle $\FF$ is {\bf counital} if $(\id\otimes_R\epsilon)\FF = 1 = (\epsilon\otimes_R\id)\FF$. Denoting $\FF =
\sum_i \FF^{1i}\otimes\FF^{2i}:=\FF^1\otimes\FF^2$, the counitality can be rewritten as $\beta(\epsilon(\FF^2))\FF^1 = 1 = \alpha(\epsilon(\FF^1))\FF^2$ and the 2-cocycle identity as
\begin{equation}\label{eq:2coccomp}
  F^1_{(1)}F^{1'}\otimes_R F^1_{(2)}F^{2'}\otimes_R F^{2} =
  F^1 \otimes_R F^2_{(1)}F^{1'} \otimes F^2_{(2)}F^{2'}
\end{equation}
where we use prime indices to distinguish different summation pairs on the same side.

We say that $\FF\in H\otimes_R H$ is {\bf invertible} if there exists an element $\tilde{\bar{\FF}}\in H\otimes_k H$
such that for one (and then for every) lift
$\tilde{\FF}\in H\otimes_k H$ of $\FF$,
$$
\tilde{\FF}\tilde{\bar{\FF}}\in 1\otimes_\genfd 1 + I_R,
\,\,\,\,\,\,
\tilde{\bar{\FF}}\tilde{\FF}\in
   1\otimes_\genfd 1 + \tilde{\bar\FF} I_R.
   $$
   Now the classes $\FF^{-1} := \tilde{\bar{\FF}}+\tilde{\bar{F}} I_R$
   and $\FF^{-1}I_R = \tilde{\bar{F}}I_R$
   do not depend on the choice of $\tilde{\bar{\FF}}$ and the
   invertibility can be rewritten as
   $$
\FF\FF^{-1} = 1\otimes_\genfd 1 + I_R = 1\otimes_R 1,
\,\,\,\,\,\,
\FF^{-1}\FF = 1\otimes_\genfd 1 + \FF^{-1} I_R.
$$
This definition is different from
but equivalent to the definition of invertibility of a 2-cocycle by Xu.
For $R = k$ bialgebroids are $k$-bialgebras and the invertible counital 2-cocycles are also known as Drinfeld twists. In general case, the invertible counital 2-cocycles are also called Drinfeld--Xu 2-cocycles or twists. 

\nxpoint As shown by Xu~\cite{xu} an $R$-bialgebroid $(H,\alpha,\beta,\Delta,\epsilon)$ with an invertible counital 2-cocycle $\FF$ induces an $R_{\FF}$-bialgebroid $(H,\alpha_{\FF},\beta_{\FF},\Delta_{\FF},\epsilon_{\FF})$ (also called the $\FF$-twist) such that $R_{\FF} = R$ as a $k$-module, equipped with an associative multiplication $\star_{\FF}$ given by $f\star_{\FF} g = \mu_R \FF(f\otimes g)$. Algebra $H = (H,\mu_H)$ does not change, the new source and target maps are given by
\begin{equation}
  \alpha_{\FF}(r) = \alpha(F^{1}\triangleright r)F^{2}
\end{equation}
\begin{equation}
  \beta_{\FF}(r) = \beta(F^{2}\triangleright r)F^{1}
\end{equation}
and the two provide an $R_{\FF}$-bimodule structure on $(H,\mu)$. The kernel $I_{R_\FF} :=\operatorname{Ker}(H\otimes_k H\to H\otimes_{R_{\FF}}H)$ of the projection is a left ideal in $H\otimes_\genfd H$ generated by $\{\beta_\FF(r)\otimes_\genfd 1-1\otimes_\genfd\alpha_\FF(r)\,|\, r\in R\}$ and $I_{R_\FF}=\FF^{-1}I_R$. 
Then the inverse $\FF^{-1}\in H\otimes_{R_{\FF}}H$ and the Sweedler like notation $\FF^{-1} = \bar{F}^{1}\otimes_{R_{\FF}}\bar{F}^{2}:= \sum_i  \bar{F}^{1i}\otimes_{R_{\FF}}\bar{F}^{2i}$ will be used.

The new comultiplication $\Delta_{\FF}:H\to H\otimes_{R_{\FF}}H$ is given by
$$
H\ni h \mapsto \Delta_{\FF} (h)  = \FF^{-1}\Delta(h)\FF, 
$$
with the counit $\epsilon_\FF$ which is the same $\genfd$-module map $\epsilon$ reinterpreted as an $R_{\FF}$-bimodule map $\epsilon_{\FF}:H\to R_{\FF}$.

The following observations will be freely used in the paper.

 \nxpoint {\bf Observation.}
 Any identity of the form $\sum_{i=1}^k h_{i} \otimes_R g_{i} = 0$ is
equivalent to the identity $\sum_{i=1}^k (\bar{F}^{1}h_{(i)})\otimes_{R_{\FF}}(\bar{F}^{2}g_{(i)}) = 0.$

Conversely, $\sum_{j=1}^l h_{j}'\otimes_{R_{\FF}}g_{j}' = 0$ is equivalent to
$\sum_{j=1}^l F^{1}h_{j}\otimes_R F^{2}g_{j} = 0$. 

\nxpoint {\bf Observation.} \label{x1hytx2gz}
$H\otimes_R H$ is not an algebra in general. However, if 
$\sum_{i=1}^k h_{i} \otimes_R g_{i} = 0$ and $x,y,z\in H$, then
the fact that $I_R$ is a right ideal and $\Delta(x)I_R\subset I_R$,
imply that the following sum is well-defined and vanishes:
$$
\sum_{i=1}^k x_{(1)} h_i y\otimes_R x_{(2)} g_i z = 0.
$$
\section{Invertibility of $(S\FF^1)\FF^2$ and map $S_{\FF}$}

\nxpoint In the situation above,
\begin{equation}
  V_{\FF} := (S\FF^{1})\FF^{2}
\end{equation}
is a well-defined element in $H$. If $H$ is a Hopf algebra then a standard calculation shows that $\overline{\FF}^{1}(S\overline{\FF}^{2})$ is the two-sided inverse of $V$ with respect to the multiplication in $H$. The calculation does not extend to Hopf algebroids, namely not only that  $\overline{\FF}^{1}(S\overline{\FF}^{2})$ is not the inverse of $V$, but worse, $\overline{\FF}^{1}(S\overline{\FF}^{2})$ is even not a well-defined expression because $\mu(\id\otimes S)I_{R_\FF}\neq 0$. We do not know if the inverse of $V_{\FF}$ exists in general.

\nxpoint Let $H$ be a left $R$-bialgebroid
and $S:H\to H$ any algebra antihomomorphism
such that $S\circ\beta = \alpha$ and $V_{\FF} = (S\FF^{1})\FF^{2}$
has an inverse $V^{-1}_{\FF}$ in $H$. Then define
\begin{equation}\label{eq:SF}
S_{\FF} h = V^{-1}_{\FF} (S h) V_{\FF} = V^{-1}_{\FF} (S h) (S\FF^{1})\FF^{2}
\end{equation}
Formula $h\mapsto S_{\FF} h$ then defines an antihomomorphism of algebras
$S_{\FF}:H\to H$ and
\begin{equation}\label{eq:VmF}
V^{-1}_{\FF}= (S_{\FF}\overline{\FF}^{1})\overline{\FF}^{2},
\end{equation}
where the right hand side is well-defined and in particular 
$\mu(S_{\FF}\otimes\id) I_{R_\FF} = 0$.

$S$ is an antihomomorphism hence it is clear that $S_{\FF}$ given by~(\ref{eq:SF}) is an antihomomorphism as well.
Regarding that $\mu(S\otimes\id)I_R = 0$ and
$F^{1'}\bar{F}^{1}\otimes_R\bar{F}^{2'}F^{2}= \FF\FF^{-1} = 1\otimes_k 1 + I_R$ we obtain 
$(S\overline{\FF}^{1})(S\FF^{1})\FF^{2}\overline{\FF}^{2} = 1$. The primed Sweedler indices $1',2'$ now refer to another copy of $\FF$ (if we have $\bar{F}$-s and $F$-s they are viewed automatically as belonging to different copies).
$$\begin{array}{lcl}
(S_{\FF}\overline{\FF}^{1})\overline{\FF}^{2}(S\FF^{1})\FF^{2} & = & 
V^{-1}_{\FF}(S\overline{\FF}^{1})V_{\FF}\overline{\FF}^{2}(S\FF^{1})\FF^{2}
\\  & = & V^{-1}_{\FF}(S\overline{\FF}^{1})(S\FF^{1'})\FF^{2'}\overline{\FF}^{2}(S\FF^{1})\FF^{2}
\\  & = & V^{-1}_{\FF}(S\FF^{1})\FF^{2}
\\  & = & 1
\end{array}$$

$$\begin{array}{lcl}
(S\FF^{1})\FF^{2}(S_{\FF}\overline{\FF}^{1})\overline{\FF}^{2} & = & (S\FF^{1})\FF^{2}V^{-1}_{\FF}(S\overline{\FF}^{1})V_{\FF}\overline{\FF}^{2}
\\  & = & (S\FF^{1})\FF^{2}V^{-1}_{\FF}(S\overline{\FF}^{1})(S\FF^{1'})\FF^{2'}\overline{\FF}^{2}
\\  & = & (S\FF^{1})\FF^{2}V^{-1}_{\FF}
\\  & = & 1
\end{array}$$

Regarding that the two-sided inverse in an associative algebra is unique, we conclude that $V_{\FF}^{-1}$ with the right hand side in~(\ref{eq:VmF}) is well-defined.

\nxpoint\label{pt:SFinv}
              {\bf Proposition.} If $S:H\to H$ is an invertible
              antihomomorphism satisfying $S\circ\beta = \alpha$,
              then the map $S_\FF$ given by~(\ref{eq:SF}) is invertible
              as well, with the inverse $S_\FF^{-1}$ given by
\begin{equation}\label{eq:SFm}
S_\FF^{-1}h = (S^{-1}V_\FF^{-1}) (S^{-1}h) (S^{-1}V_\FF),\,\,\,\,h\in H,
\end{equation}
where
\begin{equation}\label{eq:SmVF}
S^{-1} V_\FF = S^{-1}((S F^{1})F^{2}) = (S^{-1}F^{2})F^{1},
\end{equation}
\begin{equation}\label{eq:SmVFm}
S^{-1}V_\FF^{-1} = (S^{-1}V_\FF)^{-1} =
(S^{-1}_\FF\overline{F}^{2})\overline{F}^{1}
\end{equation}

{\it Proof.} If $V_\FF$ is invertible and $S$ bijective antihomomorphism,
clearly $S^{-1}V_\FF$ is also invertible with inverse $S^{-1}V_\FF^{-1}$.
Now,
$$\begin{array}{lcl}
S^{-1}_\FF(S_\FF h) &=&  (S^{-1}V_\FF^{-1}) S^{-1}(V_\FF^{-1}(Sh) V_\FF) S^{-1}V_\FF
\\
&=&
(S^{-1}V_\FF^{-1}) (S^{-1}V_\FF)(S^{-1}S h)(S^{-1} V_\FF^{-1})(S^{-1}V_\FF) = h
\end{array}$$
On the other hand,
$$\begin{array}{lcl}
  S_\FF(S_\FF^{-1} h) &=& V_\FF^{-1} S(S_\FF^{-1} h) V_\FF
  \\
  &=& V_\FF^{-1} V_\FF h  V_\FF^{-1} V_\FF = h
\end{array}$$
Therefore, equation~(\ref{eq:SFm}) follows. Knowing that the two sided inverse exists, it is unique, and for~(\ref{eq:SmVFm}) it is enough to show that $(S_\FF^{-1}\overline\FF^{2})\overline{\FF}^{1}$ is a left inverse of $S^{-1}V_\FF$. For example,
$$\begin{array}{lcl}
(S_{\FF}^{-1}\overline\FF^{2})\overline\FF^{1}(S^{-1}V_\FF) & = & 
(S^{-1}V^{-1}_\FF)(S^{-1}\overline{\FF}^{2})(S^{-1}V_\FF)\overline\FF^{1}(S^{-1}V_\FF)
\\  & = & (S^{-1}V^{-1}_\FF)(S^{-1}\overline\FF^{2})(S^{-1}\FF^{2'})\FF^{1'}\overline{\FF}^{1}(S^{-1}V_\FF)
\\  & = & (S^{-1}V^{-1}_\FF)(S^{-1}V_\FF)
\\  & = & S^{-1}(V_\FF V_\FF^{-1}) = 1
\end{array}$$
In the third step we used $(S^{-1}\overline\FF^{2})(S^{-1}\FF^{2'})\FF^{1'}\overline{\FF}^{1}=1$. To prove it one uses $\FF\FF^{-1} = 1\otimes_R 1$ and then flips the tensor factors using the flip map $H\otimes_R H\to H^\op\otimes_{R^\op}H^\op$ and applies map
$\mu^\op(\id\otimes S^{-1}):H^\op\otimes_{R^\op}H^\op\to H^\op$
which is well-defined because $S^{-1}\circ\alpha=\beta$.

\section{$S_\FF$ is an antipode map for $H_\FF$}

One would like to conclude that $S_\FF$ is an antipode for the twisted bialgebroid. The standard proofs for the Hopf algebras do not seem to generalize in straightforward manner. We first need to make some technical preparations to calculate within the new context. 

\nxpoint {\bf Lemma.}\label{pt:Sx1ytSx2z}
There are unique and well-defined $\genfd$-linear
maps $\delta_S,\delta_{S^{-1}}:H\otimes_R H\otimes_R H\to H\otimes_R H$ such that for all $x,y,z\in H$, 
\begin{equation}\label{eq:deltaS}
  \delta_S :x\otimes_R y\otimes_R z\mapsto (S x)_{(1)}y\otimes_R (S x)_{(2)}z
\end{equation}
\begin{equation}\label{eq:deltaSm}
  \delta_{S^{-1}}:x\otimes_R y\otimes_R z\mapsto (S^{-1} z)_{(1)}x\otimes_R(S^{-1}z)_{(2)}y
\end{equation}

{\it Proof.} For $\delta_S$, we start with the map $H\otimes_\genfd H\otimes_R H\to H\otimes_R H$ defined by $x\otimes_\genfd y\otimes_\genfd z\mapsto (S x)_{(1)}y\otimes_R (S x)_{(2)}z = \Delta(S x)\cdot(y\otimes_R z)$. Regarding that the elements in the image of $\Delta$ preserve the ideal $I_R$,
$\Delta(H) I_R \subset I_R$ this is well-defined. Now we need also to show that this map vanishes on the elements of the form $\beta(a)x\otimes y\otimes_R z - x\otimes\alpha(a)y\otimes_R z$. This follows from the calculation
$$\begin{array}{lcl}
  \Delta(S(\beta(a)x))& =& \Delta(S(x)\alpha(a)) =
\Delta(S(x))(\alpha(a)\otimes_R 1)\\
&=& (S x)_{(1)}\alpha(a)\otimes_R(S x)_{(2)},\\
\Delta(S(\beta(a)x))\cdot(y\otimes_R z) &=& \Delta(S(x))\cdot(\alpha(a)y\otimes_R z).
\end{array}$$
Similarly for $\delta_{S^{-1}}$ we start with the map $H\otimes_R H\otimes_\genfd H\to H\otimes_R H$ defined by $x\otimes_\genfd y\otimes_\genfd z\mapsto (S^{-1} z)_{(1)}x\otimes_R (S^{-1} z)_{(2)}y = \Delta(S^{-1}z)\cdot(x\otimes_R y)$ and show that it vanishes on the elements of the form  $x\otimes_R \beta(b) y\otimes z - x\otimes_R y\otimes \alpha(b) z$ by showing 
$$
\Delta(S^{-1}(\alpha(b)z))(x\otimes_R y) = \Delta(S^{-1}z)(\beta(b)x\otimes_R y).
$$

\nxsubpoint Maps $\delta_S,\delta_{S^{-1}}$ will be used multiple times in the proofs, in an essential way. Notice also that~(\ref{eq:antShS}) may be restated as $\delta_S(\Delta(h)\otimes_R 1_H) = 1_H\otimes_R S h$ and (\ref{eq:antSmSmh}) as $\delta_{S^{-1}}(1_H\otimes_R \Delta(h)) = S^{-1} h\otimes_R 1$. From these we
easily obtain identities
\begin{equation}\label{eq:SF11F12t}
(S F^1_{(1)})_{(1)} F^1_{(2)}\otimes_R(S F^1_{(1)})_{(2)} F^2 = 1\otimes_R (S F^1)F^2
\end{equation}
\begin{equation}\label{eq:SmF22F1}
(S^{-1}F^2_{(2)})_{(1)} F^1\otimes_R (S^{-1}F^2_{(2)})_{(2)} F^2_{(1)} = (S^{-1}F^2)F^1\otimes_R 1_H.
\end{equation}

\nxpoint \label{alphaSlemma}
    {\bf Lemma.} {\em The following identity holds for all $b\in R$ and
  $h\in H$,
      \begin{equation}\label{eq:alphaSlemma}
        \alpha((S h_{(1)})_{(1)} h_{(2)}\triangleright b) (S h_{(1)})_{(2)}
        = \alpha(b) S h.
        \end{equation}
  }

    {\it Proof.} Equation~(\ref{eq:alphaSlemma}) follows by applying the map $x\otimes_R y \mapsto \alpha(x\triangleright b)\cdot y$ to the identity~(\ref{eq:antShS}) provided this map is well-defined. In other words, we need to check that
    the map
    $$H\otimes_\genfd H\ni x\otimes y \mapsto \alpha(x\triangleright b) y\in H$$
factors through the projection $H\otimes_\genfd H\to H\otimes_R H$.
It is sufficient to show that for all $c\in R$, $h,h'\in H$, 
$\alpha(x\triangleright b) y = 0$
where $x\otimes y=\beta(c)h\otimes h'-h\otimes \alpha(c)h'\in I_R$.
But   
$\alpha(\beta(c)h\triangleright b) h' - \alpha(h\triangleright b)\alpha(c)h' = 0$ is immediate by~(\ref{eq:htrrbilin}).

\nxpoint {\bf Theorem.}\label{thm:SFFab}
$S_\FF\circ \beta_\FF = \alpha_\FF$.

In other words, $(S\beta_\FF(r))V_\FF = V_\FF\alpha_\FF(r)$ for all $r\in R$.

{\it Proof.} The right hand side is $V_\FF\alpha(F^1\triangleright r)F^2$. By $S\beta=\alpha$, the left hand side is
$$
S(\beta(F^2\triangleright r)F^1)V_\FF = (S F^1)\alpha(F^2\triangleright r) V_\FF.
$$
Thus, we need to prove
\begin{equation}\label{eq:stbalpha}
(S F^1)\alpha(F^2\triangleright r) (S F'^1)F'^2 = (S F'^1) F'^2\alpha(F^1\triangleright r) F^2.
\end{equation}
Starting from the right-hand side, we calculate 
$$\begin{array}{ccl}
  (S F'^1) F'^2\alpha(F^1\triangleright r) F^2
  &\stackrel{(\ref{eq:halphar})}=&
  (S F'^1)\alpha(F'^2_{(1)} F^1\triangleright r) F'^2_{(2)} F^2
  \\
  &\stackrel{(\ref{eq:2coccomp})}=&
  S (F'^1_{(1)} F^1)\alpha(F'^1_{(2)} F^2\triangleright r) F'^2
  \\
  & = & (S F^1) (S F'^1_{(1)})\alpha(F'^1_{(2)} F^2\triangleright r) F'^2 \\
  &\stackrel{(\ref{eq:halphar})}=&
(S F^1) \alpha((SF'^1_{(1)})_{(1)} F'^1_{(2)} \triangleright(F^2\triangleright r))
(SF'^1_{(1)})_{(2)} F'^2
\\
&\stackrel{(\ref{eq:alphaSlemma})}=&
(SF^1)\alpha(F^2\triangleright r) (SF'^1) F'^2.
\end{array}$$
Usage of the cocycle identity~(\ref{eq:2coccomp}) in the 2nd line,
holding only in $H\otimes_R H\otimes_R H$,
is justified by 
$\alpha(\alpha(a)\beta(b)h\triangleright r)=\alpha(a)\alpha(h\triangleright r)\alpha(b)$ (which follows from~(\ref{eq:htrrbilin}))
and $S\beta(a)=\alpha(a)$.

\nxpoint {\bf Corollary.} The expression
\begin{equation}\label{eq:lemmatwoidleft}
(S\bar{F}^1)_{(1)}(V_\FF)_{(1)}F^1\bar{F}^2\otimes_R
  (S\bar{F}^1)_{(2)}(V_\FF)_{(2)}F^2
\end{equation}
is well-defined, in other words it does not depend on the choice of a representative for $\FF^{-1}$.

{\it Proof.} By changing the representative by 
$\beta_\FF(r)h\otimes_\genfd k
- h\otimes_\genfd\alpha_\FF(r)k\in F^{-1}I_R$ we obtain the change of the expression by $A-B$ where
\begin{equation}\label{eq:Ar}
\begin{array}{lcl}
  A &=& (S h)_{(1)}(S\beta_\FF(r))_{(1)}V_{\FF(1)}F^1 k\otimes_R(S h)_{(2)}(S\beta_\FF(r))_{(2)}V_{\FF(2)}F^2 \\
  &=& \Delta((S h)(S\beta_\FF(r))V_\FF) (F^1 k\otimes_R F_2)
\end{array}
\end{equation}
and 
\begin{equation}\label{eq:Br}
\begin{array}{lcl}
  B &=& (S h)_{(1)}V_{\FF(1)}F^1\alpha_\FF(r) k\otimes_R (S h)_{(2)}V_{\FF(2)}F^2\\
  &=& (S h)_{(1)}V_{\FF(1)}F^1\alpha(F'^1\triangleright r)F'^2 k\otimes_R (S h)_{(2)}V_{\FF(2)}F^2\\
  &=& \Delta((S h) V_\FF)\alpha(F^1_{(1)}F'^1\triangleright r)F^1_{(2)}F'^2 k\otimes_R F^2\\
  &=& \Delta((S h) V_\FF)\alpha(F^1\triangleright r)F^2_{(1)}F'^1 k\otimes_R F^2_{(2)}F'^2\\
  &=& \Delta((S h) V_\FF\alpha_\FF(r))(F^1 k\otimes_R F^2)
\end{array}\end{equation}
The theorem says that $(S\beta_\FF(r))V_\FF = V_\FF\alpha_\FF(r)$ hence $A = B$.

\nxpoint {\bf Proposition.}
For any $G = G^1\otimes_{R_\FF} G^2\in H\otimes_{R_\FF}H = (H\otimes_\genfd H)/(\FF^{-1}I_R)$ the expression
$$
(S^{-1}G^2)_{(1)}(S^{-1}V_\FF)_{(1)}F^1\otimes_R (S^{-1}G^2)_{(2)}(S^{-1}V_\FF)_{(2)}F^2G^1
$$
is well-defined. In particular, 
\begin{equation}\label{eq:SmFmwd}
(S^{-1}\bar{F}^2)_{(1)}(S^{-1}V_\FF)_{(1)}F^1\otimes_R (S^{-1}\bar{F}^2)_{(2)}(S^{-1}V_\FF)_{(2)}F^2\bar{F}^1
\end{equation}
is well-defined.

{\it Proof.} 
If we change the representative of $\FF^{-1}$ by $\beta_\FF(r)h\otimes_\genfd k
- h\otimes_\genfd\alpha_\FF(r)k$,
then the expression~(\ref{eq:SmFmwd}) changes by
$\Delta(S^{-1}h)(C(r)-D(r))(1\otimes k)$ where
$$\begin{array}{l}
C(r) = (S^{-1}V_\FF)_{(1)}F^1\otimes_R (S^{-1}V_\FF)_{(2)}F^2\beta_\FF(r)h\\
D(r) = (S^{-1}\alpha_\FF(r))_{(1)}(S^{-1}V_\FF)_{(1)} F^1\otimes_R
(S^{-1}\alpha_\FF(r))_{(2)}(S^{-1}V_\FF)_{(2)} F^2.
\end{array}$$
If we apply antihomomorphism $S^{-1}$ to 
$V_\FF\alpha_\FF(r)=(S\beta_\FF(r))V_\FF$ we obtain
\begin{equation}
(S^{-1}\alpha_F)(r)S^{-1}V_\FF = (S^{-1}V_\FF)\beta_\FF(r).
\end{equation}
Thus,
$$\begin{array}{lcl}
  D(r)&=& \Delta((S^{-1}V_\FF)\beta(F^2\triangleright r)F^1)\FF\\
  &=& \Delta(S^{-1}V_\FF)(F^1_{(1)}F'^1\otimes_R\beta(F^2\triangleright r)
  F^1_{(2)} F'^2)\\
  &=& \Delta(S^{-1}V_\FF)(F^1\otimes_R\beta(F^2_{(2)}F'^2\triangleright r)F^2_{(1)} F'^1)\\
  &=& \Delta(S^{-1}V_\FF)(F^1\otimes F^2\beta(F'^2\triangleright r)F'^1)\\
  &=& \Delta(S^{-1}V_\FF)(F^1\otimes F^2\beta_\FF(r))\\
  &=& C(r).
\end{array}$$

\nxpoint {\bf Lemma.} The following identities hold in $H\otimes_R H$,
\begin{equation}\label{eq:lemmatwoid}
(S\bar{F}^1)_{(1)}(V_\FF)_{(1)}F^1\bar{F}^2\otimes_R
  (S\bar{F}^1)_{(2)}(V_\FF)_{(2)}F^2 = 1\otimes_R V_\FF
\end{equation}
\begin{equation}\label{eq:lemmatwoid2}
(S^{-1}\bar{F}^2)_{(1)}(S^{-1}V_\FF)_{(1)}F''^1\otimes_R (S^{-1}\bar{F}^2)_{(2)}(S^{-1}V_\FF)_{(2)}F''^2\bar{F}^1 = S^{-1}V_\FF\otimes_R 1.
\end{equation}
{\it Proof.}
Applying $\delta_S$ (see~(\ref{eq:deltaS}) to the 2-cocycle condition 
$F'^1\otimes_R F'^2_{(1)}F^1\otimes_R F'^2_{(2)}F^2 
= F'^1_{(1)}F^1\otimes_R F'^1_{(2)}F^2\otimes_R F'^2$ we obtain
\begin{eqnarray*}
(SF'^1)_{(1)}F'^2_{(1)}F^1\otimes_R (SF'^1)_{(2)}F'^2_{(2)}F^2
= S(F'^1_{(1)}F^1)_{(1)}F'^1_{(2)}F^2 \otimes_R S(F'^1_{(1)}F^1)_{(2)}F'^2
\end{eqnarray*}

By~\refpoint{x1hytx2gz} if $\sum_i h_i\otimes h'_i \in I$ then 
$$
\sum_i (S\bar{F}^1)_{(1)} h_i \bar{F}^2 \otimes_R (S\bar{F}^1)_{(2)} h'_i = 0.
$$
Therefore, the left hand side in~(\ref{eq:lemmatwoid}),
temporarily denoted by $L_1$, becomes 
$$\begin{array}{lcl}
L_1 &=& (S\bar{F}^1)_{(1)}(SF'^1)_{(1)}(F'^2)_{(1)}F^1\bar{F}^2\otimes_R
(S\bar{F}^1)_{(2)}(SF'^1)_{(2)}(F'^2)_{(2)}F^2\\
&=& (S\bar{F}^1)_{(1)}(S F^1)_{(1)}S(F'^1_{(1)})_{(1)} F'^1_{(2)}
F^2\bar{F}^2 \otimes_R (S\bar{F}^1)_{(2)} (SF^1)_{(2)}(SF'^1_{(1)})_{(2)}F'^2
\end{array}$$
\begin{equation}\label{eq:SF12}
L_1 = (S (F^1\bar{F}^1))_{(1)}S(F'^1_{(1)})_{(1)} F'^1_{(2)}
(F^2\bar{F}^2) \otimes_R (S(F^1\bar{F^1}))_{(2)} (SF'^1_{(1)})_{(2)}F'^2
\end{equation}
Axiom~(\ref{eq:antShS}) for $S$ implies the identity
$$
(S F'^1_{(1)})_{(1)} F'^1_{(2)} \otimes_R (SF'^1_{(1)})_{(2)}F'^2
= 1\otimes_R (S F'^1) F'^2.
$$
Since $I_R$ is a left ideal and 
$\Delta(x)I_R\subset I_R$ a more general identity
$$
x_{(1)} (S F'^1_{(1)})_{(1)} F'^1_{(2)} y\otimes_R x_{(2)} (SF'^1_{(1)})_{(2)}F'^2
= x_{(1)}y\otimes_R x_{(2)}(S F'^1) F'^2
$$
holds for all $x,y\in H$. By linearity in $x$ and $y$, this,
along with (\ref{eq:SF12}),  gives
$$
L_1 = (S (F^1\bar{F}^1))_{(1)} (F^2\bar{F}^2)
\otimes_R (S (F^1\bar{F}^1))_{(2)} (S F'^1)F'^2.
$$
We now use $F^1\bar{F}^1\otimes_R F^2\bar{F}^2 = 1\otimes 1 + I_R$ and immediately obtain
$$
L_1 = 1\otimes_R (S F'^1)F'^2 = 1\otimes_R V_\FF
$$
hence the statement of the lemma.

\nxpoint {\bf Theorem.} \label{thm:SFF}
\begin{equation}\label{eq:SFFid}
(S_\FF h_{(1\FF)})_{(1\FF)} h_{(2\FF)}\otimes_{R_\FF} (S_\FF h_{(1\FF)})_{(2\FF)}
= 1\otimes_{R_\FF} S_\FF h
\end{equation}

{\it Proof.} To reduce~(\ref{eq:SFFid}) to something manageable, we first expand it, 
$$
(V_{\FF}^{-1})_{(1\FF)}(S h_{(1\FF)})_{(1\FF)}(V_{\FF})_{(1\FF)} h_{(2\FF)}\otimes_{R_\FF}(V_{\FF}^{-1})_{(2\FF)}(S h_{(1\FF)})_{(2\FF)}(V_{\FF})_{(2\FF)} = 1\otimes_{R_\FF} V_\FF^{-1} (S h) V_\FF
$$
Expressing $\Delta_\FF$ in terms of $\Delta$ we obtain that LHS is
$$
\bar{F'}^1(V_\FF^{-1})_{(1)}(S h_{(1\FF)})_{(1)}(V_\FF)_{(1)}F''^1 h_{(2\FF)}\otimes \bar{F'}^2(V_\FF^{-1})_{(2)} (S h_{(1\FF)})_{(2)} (V_\FF)_{(2)}F''^2
$$
We would like to check this equality in terms of the original antipode and coproduct and it is convenient to multiply the required identity with $\Delta(V_\FF)\FF$ from the left, what also changes the tensor product to $\otimes_R$ because $F I_\FF = I_R$. This changes the required identity to $L_2 = R_2$ where
\begin{equation}\label{eq:L2initial}
    L_2 = (S(\bar{F}^1 h_{(1)} F^1))_{(1)}(V_\FF)_{(1)}F''^1
    \bar{F}^2 h_{(2)} F^2\otimes_R (S(\bar{F}^1 h_{(1)} F^1))_{(2)}(V_\FF)_{(2)}F''^2 
\end{equation}
  and
$$\begin{array}{lcl}
 R_2 &=& (V_\FF)_{(1)} F^1 \otimes_R (V_\FF)_{(2)} F^2 V^{-1}_\FF (S h) V_\FF
    \\
&=& (S F'^{1})_{(1)} F'^2_{(1)} F^1\otimes_R (S F'^{1})_{(2)} F'^2_{(2)} F^2 V_\FF^{-1} (S h) V_\FF 
    \\
    &=& \delta_S(F'^{1}\otimes_R F'^2_{(1)} F^1\otimes_R F'^2_{(2)} F^2 V_\FF^{-1} (S h) V_\FF)\\
    &\stackrel{(\ref{eq:2coc})}=& \delta_S(F^1_{(1)}F'^1\otimes_R F^1_{(2)}F'^2\otimes F^2 V_\FF^{-1} (S h) V_\FF)\\
    & = & (S(F^1_{(1)}F'^1))_{(1)}F^1_{(2)}F'^2\otimes_R(S(F^1_{(1)}F'^1))_{(2)}F^2 V_\FF^{-1} (S h) V_\FF\\
    &=& \delta_S(F'^1\otimes_R (S F^1_{(1)})_{(1)} F^1_{(2)}F'^2\otimes_R(S F^1_{(1)})_{(2)}F^2 V_\FF^{-1} (S h) V_\FF) \\
    &\stackrel{(\ref{eq:SF11F12t})}=&
    \delta_S(F'^1\otimes_R F'^2\otimes_R (S F^1)F^2  V_\FF^{-1}(S h) V_\FF)\\
    &=& \delta_S(F'^1\otimes_R F'^2\otimes (S h) V_\FF)\\
    &=&  (S F'^1)_{(1)} F'^2\otimes_R  (S F'^1)_{(2)}(S h) V_\FF.
\end{array}$$

We rewrite the expression~(\ref{eq:L2initial}) for $L_2$
by using the fact that
$S$ is an antihomomorphism of corings~\refpoint{antihom}
$$\begin{array}{lcl}
    L_2  &=&
    (S F^1)_{(1)} (S h_{(1)})_{(1)}  S (\bar{F}^1_{(1)})_{(1)} (V_\FF)_{(1)} F''^1
    \bar{F}^2 h_{(2)} F^2\otimes_R (S F^1)_{(2)}(S h_{(1)})_{(2)}S(\bar{F}^1_{(2)})(V_\FF)_{(2)}F''^2
\end{array}$$
Lemma~\refpoint{pt:L2lastlemma} below and
Observation~\refpoint{x1hytx2gz} imply that 
$$\begin{array}{lcl}
  L_2 &=& (S F^1)_{(1)} (S h_{(1)})_{(1)}h_{(2)}F^2 \otimes_R (S F^1)_{(2)}(S h_{(1)})_{(2)} V_\FF\\
  &\stackrel{(\ref{eq:SF})}=& (S F^1)_{(1)}F^2\otimes_R (S F^1)_{(2)} (S h) V_\FF
  \\ &=& R_2
\end{array}$$
as required. Therefore it remains only to prove the following lemma.
\nxpoint \label{pt:L2lastlemma} {\bf Lemma.} The following identity holds and the LHS is well-defined,
\begin{equation}\label{eq:L2lastlemma}
  (S\bar{F}^1)_{(1)}(V_\FF)_{(1)}F''^1\bar{F}^2\otimes_R(S\bar{F}^1)_{(2)}(V_\FF)_{(2)}F''^2 = 1\otimes_R V_\FF
\end{equation}
{\it Proof.} Starting with the right hand side of~(\ref{eq:L2lastlemma}), we calculate
$$\begin{array}{rl}
  (S\bar{F}^1)_{(1)}& \!\!\!\!\!(V_\FF)_{(1)} F''^1\bar{F}^2
  \otimes_R (S\bar{F}^1)_{(2)}(V_\FF)_{(2)}F''^2
  =
  \\ &=
(S\bar{F}^1)_{(1)}(S F'^1)_{(1)}F'^2_{(1)}F''^1\bar{F}^2\otimes_R (S\bar{F}^1)_{(2)}(S F'^1)_{(2)}F'^2_{(2)}F''^2
  \\
  &= \delta_S\left(\bar{F}^1\otimes_R(S F'^1)_{(1)}F'^2_{(1)}F''^1\bar{F}^2\otimes_R (S F'^1)_{(2)}F'^2_{(2)}F''^2\right)
  \\
  &= \delta_S\left(\bar{F}^1\otimes_R\delta_S(F'^1\otimes_R F'^2_{(1)}F''^1\bar{F}^2\otimes_R F'^2_{(2)}F''^2)\right)
  \\
  &= \delta_S\left(\bar{F}^1\otimes_R\delta_S(F'^1_{(1)}F''^1\otimes_R
  F'^1_{(2)}F''^2\bar{F}^2\otimes_R F'^2)\right)
  \\
  &=
  (S\bar{F}^1)_{(1)}(S (F'^1_{(1)}F''^1))_{(1)} F'^1_{(2)}F''^2\bar{F}^2\otimes_R (S\bar{F}^1)_{(2)}(S(F'^1_{(1)}F''^1))_{(2)} F'^2
  \\ &=
  (S(F'^1_{(1)}F''^1\bar{F}^1)_{(1)}F'^1_{(2)}(F''^2\bar{F}^2)\otimes_R(S(F'^1_{(1)}F''^1\bar{F}^1))_{(2)}F'^2 \\
  &=\delta_S(F'^1_{(1)}F''^1\bar{F}^1\otimes_R F'^1_{(2)}(F''^2\bar{F}^2)\otimes_R F'^2)\\
  &=\delta_S(F'^1_{(1)}\otimes_R F'^1_{(2)}\otimes_R F'^2)\\
  &= (S F^1_{(1)})_{(1)}F^1_{(2)}\otimes_R(S F^1_{(1)})_{(2)}F^2
  \\
  &\stackrel{(\ref{eq:SF11F12t})}= 1\otimes_R (S F^1)F^2
  \\
  &= 1\otimes_R V_\FF.
\end{array}$$
\nxpoint {\bf Theorem.}\label{thm:SmFF}
\begin{equation}\label{eq:SmFFid}
(S^{-1}_\FF h_{(2\FF)})_{(1\FF)}\otimes_{R_\FF}(S^{-1}_\FF h_{(2\FF)})_{(2\FF)}h_{(1\FF)} = S^{-1}_\FF h\otimes_{R_\FF} 1_H
\end{equation}
   {\it Proof.} We multiply this equation by $\Delta(S^{-1}_\FF)\FF$ to obtain the equivalent statement $L_3 = R_3$ where
$$\begin{array}{lcl}
      L_3 &=& \Delta(S^{-1} V_\FF)\FF\left( (S^{-1}_\FF h_{(2\FF)})_{(1\FF)}\otimes_{R_\FF}(S^{-1}_\FF h_{(2\FF)})_{(2\FF)}h_{(1\FF)}\right) \\
      &=& (S^{-1} V_\FF)_{(1)}(S^{-1}_\FF h_{(2\FF)})_{(1)}F''^1\otimes_R
      (S^{-1}V_\FF)_{(2)}(S^{-1}_\FF h_{(2\FF)})_{(2)}F''^2 h_{(1\FF)}
      \\
 &\stackrel{(\ref{eq:SFm})}=& (S^{-1}h_{(2\FF)})_{(1)}(S^{-1} V_\FF)_{(1)}F''^1\otimes_R
      (S^{-1} h_{(2\FF)})_{(2)}(S^{-1}V_\FF)_{(2)}F''^2h_{(1\FF)}
\end{array}$$
and
$$
      R_3 =  \Delta(S^{-1} V_\FF)\FF(S^{-1}_\FF h\otimes_{R_\FF} 1_H)\\
$$
Now, recall~(\ref{eq:SmVF}) that $S^{-1}V_\FF = (S^{-1}F^2)F^1$, hence
$$\begin{array}{lcl}
  R_3 &=& (S^{-1}F^2)_{(1)}F^1_{(1)}F'^1(S^{-1}_\FF h)
  \otimes_R (S^{-1}F^2)_{(2)}F^1_{(2)}F'^2\\
  &\stackrel{(\ref{eq:deltaSm})}=&
  \delta_{S^{-1}}\left(F^1_{(1)}F'^1(S^{-1}_\FF h)\otimes_R F^1_{(2)}F'^2\otimes_R F^2\right)\\
  &\stackrel{(\ref{eq:2coccomp})}=&
\delta_{S^{-1}}(F^1(S^{-1}_\FF h)\otimes_R F^2_{(1)}F'^1\otimes_R F^2_{(2)}F'^2)
\\
&=& (S^{-1}F'^2)_{(1)}(S^{-1}F^2_{(2)})_{(1)} F^1(S^{-1}_\FF h)
\otimes_R(S^{-1}F'^2)_{(2)}(S^{-1}F^2_{(2)})_{(2)} F^2_{(1)} F'^1
\end{array}$$
Using $(S^{-1}F^2)_{(1)} F^1\otimes_R (S^{-1}F^2)_{(2)}F^2_{(1)} \stackrel{(\ref{eq:SmF22F1})}= (S^{-1}F^2)F^1\otimes_R 1 = S^{-1}V_\FF\otimes_R 1$ we conclude
$$
R_3 = (S^{-1}F'^2)_{(1)}(S^{-1}V_\FF)(S^{-1}_\FF h)\otimes_R(S^{-1}F'^2)_{(2)}F'^1,
$$
hence by using $(S^{-1}V_\FF)(S^{-1}_\FF h) = (S^{-1}h)(S^{-1}V_\FF)$,
\begin{equation}\label{eq:R3final}
R_3 = (S^{-1}F'^2)_{(1)}(S^{-1}h)(S^{-1}V_\FF)\otimes_R(S^{-1}F'^2)_{(2)}F'^1. 
\end{equation}
Now we perform the calculation for $L_3$.
$$\begin{array}{lcl}
  L_3 &=& (S^{-1}(\bar{F}^2 h_{(2)}F^2))_{(1)}(S^{-1}V_\FF)_{(1)}F'^1\otimes_R
  (S^{-1}(\bar{F}^2 h_{(2)}F^2))_{(2)}(S^{-1}V_\FF)_{(2)}F'^2\bar{F}^1 h_{(1)} F^1
  \\
  &=&
  (S^{-1}(h_{(2)}F^2))_{(1)}(S^{-1}\bar{F}^2)_{(1)}(S^{-1}V_\FF)_{(1)}F'^1
  \otimes_R\\
  && \otimes_R
  (S^{-1}(h_{(2)}F^2))_{(2)}(S^{-1}\bar{F}^2)_{(2)}(S^{-1}V_\FF)_{(2)}F'^2\bar{F}^1 h_{(1)}F^1\\
  &=&\delta_{S^{-1}}\left((S^{-1}\bar{F}^2)_{(1)}(S^{-1}V_\FF)_{(1)}F'^1\otimes_R
   (S^{-1}\bar{F}^2)_{(2)}(S^{-1}V_\FF)_{(2)}F'^2 \bar{F}^1 h_{(1)}F^{1}
  \otimes_R h_{(2)} F^2\right)\\
  &\stackrel{\refpoint{pt:lastlemma}}=&
  \delta_{S^{-1}}\left(S^{-1}V_\FF\otimes_R h_{(1)}F^1\otimes_R h_{(2)}F^2\right)
  \\
  &=&(S^{-1}F^2)_{(1)}(S^{-1}h_{(2)})_{(1)}S^{-1}V_\FF\otimes_R (S^{-1}F^2)_{(2)}(S^{-1}h_{(2)})_{(2)}h_{(1)}F^1
  \\
  &=&\delta_{S^{-1}}\left((S^{-1}h_{(2)})_{(1)}(S^{-1}V_\FF)\otimes_R(S^{-1}h_{(2)})_{(2)}h_{(1)}F^1 \otimes_R  F^2\right)\\
  &\stackrel{(\ref{eq:antSmSmh})}=&
  \delta_{S^{-1}}\left((S^{-1}h)(S^{-1}V_\FF)\otimes_R F^1\otimes_R F^2\right)\\
  &=& (S^{-1}F^2)_{(1)} (S^{-1}h)(S^{-1}V_\FF)\otimes_R (S^{-1}F^2)_{(2)} F^1\\
  &=& R_3.
\end{array}$$

The proof is complete modulo the following lemma. 
\nxpoint \label{pt:lastlemma} {\bf Lemma.} 
\begin{equation}\label{eq:Smlem}
(S^{-1}\bar{F}^2)_{(1)}(S^{-1}V_\FF)_{(1)}F'^1\otimes_R (S^{-1}\bar{F}^2)_{(2)}(S^{-1}V_\FF)_{(2)}F'^2\bar{F}^1 = S^{-1}V_\FF\otimes_R 1.
\end{equation}
{\it Proof.} We write the 2-cocycle condition in the form
$$
F^1_{(1)}F'^1\otimes_R F^1_{(2)}F'^2\otimes_R F^2 = F'^1\otimes_R F'^2_{(1)}F^1\otimes_R F'^{2}_{(2)}F^2
$$
and apply $\delta_{S^{-1}}$ to it to obtain
$$
(S^{-1}F^2)_{(1)}F^1_{(1)}F'^1\otimes_R (S^{-1}F^2)_{(2)}F^1_{(2)}F'^2
= (S^{-1}F^2)_{(1)}(S^{-1}F'^2)_{(1)}F'^1\otimes_R(S^{-1}F^2)_{(2)}(S^{-1}F'^2)_{(2)}F'^2F^1.
$$
hence after using~(\ref{eq:SmF22F1}) for the RHS,
$$
(S^{-1}V_\FF)_{(1)}F'^1\otimes_R (S^{-1}V_\FF)_{(2)}F'^2
 = (S^{-1}F^2)_{(1)}(S^{-1}V_\FF)\otimes_R(S^{-1}F^2)_{(2)}F^1
$$
 Tensor multiply over $\genfd$ this equation by $F^{-1}=\overline{F}^1\otimes_R\overline{F}^2$ from the right and apply the multiplication $H\otimes_\genfd H\to H$ of the two middle factors to obtain
 $$\begin{array}{l}
(S^{-1}V_\FF)_{(1)}F'^1\otimes_R (S^{-1}V_\FF)_{(2)}F'^2\overline{F}^1\otimes_R\overline{F}^2 =\\
   = (S^{-1}F^2)_{(1)}(S^{-1}V_\FF)\otimes_R(S^{-1}F^2)_{(2)}F^1\overline{F}^1\otimes_R\overline{F}^2,
   \end{array}
 $$
 and then apply $\delta_{S^{-1}}$ to obtain~(\ref{eq:Smlem}) after the following simplification at the right hand side: 
 $(S^{-1}\overline{F}^2)_{(1)}(S^{-1}F^2)_{(1)}\otimes_R (S^{-1}\overline{F}^2)_{(2)}(S^{-1}F^2)_{(2)}F^1\overline{F}^1=\delta_{S^{-1}}(1\otimes_R F^{1}\overline{F}^1\otimes_R F^2\overline{F}^2)=\delta_{S^{-1}}(1\otimes_R1\otimes_R 1) = 1\otimes_R 1$.

\nxpoint {\bf Main theorem.} Let $H = (H,\alpha,\beta,\Delta,\epsilon,S)$ be a Hopf $R$-algebroid with an invertible antipode $S$ and $\FF = \sum F^1\otimes F^2\in H\otimes_R H$ an invertible counital 2-cocycle. 
If $V_\FF = (S F^1)F^2$ is invertible then the algebra antihomomorphism 
$S_\FF:h\mapsto V_\FF^{-1}(S h) V_\FF$ is an invertible antipode
for the twisted $R_\FF$-bialgebroid
$H_\FF = (H,\alpha_\FF,\beta_\FF,\Delta_\FF,\epsilon_\FF)$ and its 
inverse is $S_\FF^{-1}:h\mapsto (S^{-1}V_\FF^{-1})(S^{-1} h)(S^{-1}V_\FF)$.

{\it Proof.} This follows from~\refpoint{pt:SFinv}, the fact that $S_\FF$ is an antihomomorphism of algebras, and theorems \refpoint{thm:SFFab}, \refpoint{thm:SFF} and \refpoint{thm:SmFF}.

\nxpoint In B\"ohm's formalism~\cite{bohmnew,bohmHbk} of a symmetric Hopf algebroid (with related left and right Hopf algebroids) an antipode is unique (just like in te Hopf algebra case), but in the formalism of a left bialgebroid and invertible antipode as above there is a (well studied) non-uniqueness of the antipode. It is well known that whenever $\FF\in H\otimes_R H$ is an invertible bialgebroid 2-cocycle, then $\FF^{-1}\in H_\FF\otimes_{R_\FF}H_\FF$ is a bialgebroid 2-cocycle. According to~(\ref{eq:VmF}), $V_{\FF^{-1}} = V_\FF^{-1}$, hence conjugation with $V_{\FF^{-1}}$ is the inverse of the operation of conjugation with $V_\FF$. Thus the twisting by the inverse twist recovers back not only the original bialgebroid but also the original (choice of) antipode; the conjugution formula gives a bijection between the sets of all possible antipodes for the original and for the twisted bialgebroids.  

\nxpoint {\it On examples.} N. Kowalzig (\cite{kowalzigth}) has shown the conditions when the universal enveloping algebra of a Lie algebroid which is a left bialgebroid in fact has an antipode. The condition is the existence of certain type of a connection. For the tangent Lie algebroid, the universal enveloping bialgebroid is that of differential operators. Ping Xu~\cite{xu} has proved that any deformation quantization of a Poisson manifold provides a twist for bialgebroid of differential operators with coefficients in formal power series in one variable.  The invertibility of $V_\FF$ is easy to check: it is a formal power series whose constant term is $1$. In one of the simplest cases, the twisted Hopf algebroid is treated explicitly in~\cite{halgoid} and we explicitly checked in~\cite{twlinPois} that the formula for the antipode proved in this article in fact gives the (known from~\cite{halgoid}) antipode for the twisted bialgebroid. That affine space example is (a completed version of) a scalar Hopf algebroid~\cite{stojicscalar} and the twist is related to a Drinfeld twist of the underlying Hopf algebra.

Drinfeld type twisting is well known for weak Hopf algebras as well. Each weak Hopf algebra induces canonically a Hopf algebroid. The twist for a weak Hopf algebra induces an invertible 2-cocycle and $V_\FF$ is also automatically invertible in that case.

It would be interesting to study the generalization of the results in this paper for the dynamical twists. It is also not clear how to define the twisting cocycles for the symmetric Hopf algebroids (in the generality when the antipode is not necessarily invertible). One expects to have a pair of bialgebroid 2-cocycles -- one for the left and one for the right bialgebroid -- with some compatibility conditions. 

\bibliographystyle{alpha}

\end{document}